\documentclass[reqno,a4paper]{amsart}
\usepackage[latin1]{inputenc}
\usepackage[T1]{fontenc}
\usepackage[pdfborder={0 0 0}]{hyperref}
\usepackage{lmodern,amssymb}

\newtheorem{thm}{Theorem}
\newtheorem{lem}[thm]{Lemma}
\newtheorem{cor}[thm]{Corollary}
\theoremstyle{remark}
\newtheorem{remark}[thm]{Remark}
\let\epsilon\varepsilon
\let\phi\varphi
\newcommand{\bbR}{\mathbb{R}}
\newcommand{\TV}{\operatorname{TV}}
\newcommand{\abs}[2][]{\mathopen#1|#2\mathclose#1|}
\newcommand{\norm}[2][]{\mathopen#1\|#2\mathclose#1\|}
\newcommand{\setof}[2][]{\mathopen#1\{#2\mathclose#1\}}
\newcommand{\suchthat}{\colon} 
\DeclareMathOperator{\dist}{dist}
\newcommand{\loc}{_{\rm loc}}

\newcommand{\cond}[1]{\textup{(\textit{#1})}}

\title{The Kolmogorov--Riesz compactness theorem}

\author[Hanche-Olsen]{Harald Hanche-Olsen}
\address[Hanche-Olsen]{\newline
    Department of Mathematical Sciences\newline
    Norwegian University of Science and Technology\newline
    NO--7491 Trondheim, Norway}
\email{hanche@math.ntnu.no}
\urladdr{http://www.math.ntnu.no/\~{}hanche/}

\author[Holden]{Helge Holden}
\address[Holden]{\newline
   Department of Mathematical Sciences\newline
   Norwegian University of Science and Technology\newline
   NO--7491 Trondheim, Norway,\newline
{\rm and} \newline
 Centre of Mathematics for Applications\newline
  University of Oslo\newline
 P.O.\ Box 1053, Blindern\newline
 NO--0316 Oslo, Norway }
\email[]{holden@math.ntnu.no}
\urladdr{http://www.math.ntnu.no/\~{}holden/}

\date{\today}
\begin{document}
\thanks{Supported in part by the Research Council of Norway. 
This paper was written as part of  the international research 
program on Nonlinear Partial Differential Equations at 
the Centre for Advanced Study at the Norwegian Academy 
of Science and Letters in Oslo during the academic year 2008--09.
The authors would like to thank Helge Kristian Jenssen for fruitful
discussions.}

  \begin{abstract}
We show that the Arzelà--Ascoli theorem and Kolmogorov compactness theorem
both are
consequences of a simple lemma on compactness in metric spaces.  Their
relation to Helly's theorem is discussed. The paper contains a
detailed discussion on the historical background of the Kolmogorov
compactness theorem.  
  \end{abstract}
\maketitle
\section{Introduction}
Compactness results in the spaces $L^p(\bbR^d)$ ($1\le p<\infty$) 
are often vital in
existence proofs for nonlinear partial differential equations.
A necessary and sufficient condition for a subset of $L^p(\bbR^d)$ to
be compact is
given in what is often called the Kolmogorov compactness theorem, or 
Fréchet--Kolmogorov compactness theorem. Proofs of this theorem are
frequently based on the Arzelà--Ascoli theorem. We here show how one can
deduce both the Kolmogorov compactness theorem and the Arzelà--Ascoli theorem
from one common lemma on compactness in metric spaces, which again is
based on the fact that a metric space is compact if and only if it is
complete and totally bounded.

Furthermore, we trace out the historical roots of Kolmogorov's
compactness theorem, which originated in Kolmogorov's 
classical paper \cite{Kolmogorov} from 1931.  However, there were
several other approaches to the issue of describing compact subsets of 
$L^p(\bbR^d)$ prior to and after Kolmogorov,  and several of these 
are described in Section~\ref{history}.
Furthermore, extensions to other spaces, say $L^p(\bbR^d)$ ($0\le
p<1$),  Orlicz spaces, or compact groups, are described.
Helly's theorem is often used as a replacement for Kolmogorov's
compactness theorem, in particular in the context of nonlinear
hyperbolic conservation laws,  in
spite of being more specialized
(e.g., in the sense that its classical version requires one
spatial dimension).
For instance, Helly's theorem is an essential ingredient in Glimm's
ground breaking existence proof for nonlinear hyperbolic systems \cite{Glimm}.
We show below that Helly's theorem is an easy
consequence of Kolmogorov's compactness theorem.

\section{Preliminary results}
An \textit{$\epsilon$-cover}
of a metric space is a cover of the space
consisting of sets of diameter at most $\epsilon$.
A metric space is called \textit{totally bounded} if it admits a
finite $\epsilon$-cover for every $\epsilon>0$.
It is well known that a metric space is compact if and only if it is
complete and totally bounded (see, e.g., \cite[p.\ 13]{Yosida}).
Since we are interested in compactness results for subsets of
Banach spaces, we may, and shall, concentrate our
attention on total boundedness.

Here is the key lemma for many compactness results
(in this lemma and its proof, every metric is named $d$):

\begin{lem} \label{lemma:key}
  Let $X$ be a metric space.
  Assume that, for every $\epsilon>0$, there exists some
  $\delta>0$, a metric space $W$, and a mapping
  $\Phi\colon X\to W$ so that $\Phi[X]$ is totally bounded, and
  whenever $x,y\in X$ are such that
  $d\big(\Phi(x),\Phi(y)\big)<\delta$, then $d(x,y)<\epsilon$.
  Then $X$ is totally bounded.
\end{lem}
\begin{proof}
  For any $\epsilon>0$, pick $\delta$, $W$ and $\Phi$ as in
  the statement of the lemma.
  Since $\Phi[X]$ is totally bounded, there exists a finite $\delta$-cover
  $\setof{V_1,\ldots,V_n}$ of $\Phi[X]$.
  Then it immediately follows from the assumptions that
  $\setof{\Phi^{-1}(V_1),\ldots,\Phi^{-1}(V_n)}$ is an $\epsilon$-cover of $X$.
  Thus $X$ is totally bounded.
\end{proof}

\noindent
Lemma \ref{lemma:key} embodies the main argument in the standard proof of
the classical Arzelà--Ascoli theorem, as we now demonstrate.

\begin{thm}[Arzelà--Ascoli] \label{thm:a}
  Let $\Omega$ be a compact topological space.
  Then a subset of $C(\Omega)$ is totally bounded
  in the supremum norm if, and only if,
  \begin{enumerate}
  \item it is pointwise bounded, and
  \item it is equicontinuous.
  \end{enumerate}
\end{thm}
\noindent
Recall the definition of equicontinuity:
Condition \cond{ii} means that
for every $x\in\Omega$ and every $\epsilon>0$
there is a neighborhood $V$ of $x$ so that
$\abs{f(y)-f(x)}<\epsilon$ for all $y\in V$
and all $f$ in the given set of functions.
\begin{proof}
  Assume $\mathcal{F}\subset C(\Omega)$ is pointwise bounded and
  equicontinuous.
  Let $\epsilon>0$.
  Combining the equicontinuity of $\mathcal{F}$ and compactness of
  $\Omega$, we can find a finite set of points
  $x_1,\ldots,x_n\in\Omega$ with neighborhoods $V_1,\ldots,V_n$
  covering all of $\Omega$ so that $\abs{f(x)-f(x_j)}<\epsilon$ whenever
  $f\in\mathcal{F}$ and $x\in V_j$.

  Define $\Phi\colon\mathcal{F}\to \bbR^n$ by
  \[\Phi(f)=\big(f(x_1),\ldots,f(x_n)\big).\]
  By the pointwise boundedness of $\mathcal{F}$, the image
  $\Phi[\mathcal{F}]$ is bounded, and hence totally bounded, in
  $\bbR^n$.

  Furthermore, if $f,g\in\mathcal{F}$ with
  $\norm{\Phi(f)-\Phi(g)}_\infty<\epsilon$, then since any
  $x\in\Omega$ belongs to some $V_j$,
  \[\abs{f(x)-g(x)}\le\abs{f(x)-f(x_j)}+\abs{f(x_j)-g(x_j)}+\abs{g(x_j)-g(x)}<3\epsilon,\]
  and so $\norm{f-g}_\infty\le3\epsilon$.
  By Lemma \ref{lemma:key}, $\mathcal{F}$ is totally bounded.

  \medskip\noindent
  For the converse, assume that $\mathcal{F}$ is a totally bounded
  subset of $C(\Omega)$.
  
  The existence of a finite $\epsilon$-cover for $\mathcal{F}$,
  for any $\epsilon$, clearly implies the boundedness of
  $\mathcal{F}$, thus establishing the uniform boundedness and hence
  also pointwise boundedness of $\mathcal{F}$.

  To prove equicontinuity, let $x\in\Omega$ and $\epsilon>0$ be given.
  Pick an
  $\epsilon$-cover $\setof{U_1,\ldots,U_n}$ of $\mathcal{F}$,
  and chose $g_j\in U_j$ for $j=1,\ldots,n$.
  Pick a neighborhood $V_j$ of $x$ so that
  $\abs{g_j(y)-g_j(x)}<\epsilon$ whenever $y\in V_j$, for $j=1,\ldots,n$.
  Let $V=V_1\cap\cdots\cap V_m$.
  If $f\in U_j$ then $\norm{f-g_j}_\infty\le\epsilon$, and so when $y\in V$,
  \[
     \abs{f(y)-f(x)}\le\abs{f(y)-g_j(y)}+\abs{g_j(y)-g_j(x)}+\abs{g_j(x)-f(x)}
     <3\epsilon,
  \]
  which proves equicontinuity.
\end{proof}

\begin{remark}
This theorem was first proved by Ascoli \cite{Ascoli} for
equi-Lipschitz functions and extended by Arzelà \cite{Arzela} to a
general family of equicontinuous functions. See \cite[p.\ 203]{Benedetto}.
\end{remark}

  We present the following theorem, first proved by Fréchet
  \cite{Frechet1} for the case $p=2$, as a warm-up exercise, as the
  proof is short and nicely exposes some key ideas for the proof of
  Theorem \ref{thm:k-r}.
\begin{thm}
  A subset of $l^p$, where $1\le p<\infty$, is totally bounded if,
  and only if,
  \begin{enumerate}
  \item it is pointwise bounded, and
  \item for every $\epsilon>0$ there is some $n$ so that, for
  every $x$ in the given subset,\[\sum_{k>n}\abs{x_k}^p<\epsilon^p.\]
  \end{enumerate}
\end{thm}
\begin{proof}
  Assume that $\mathcal{F}\subset l^p$ satisfies the two conditions.
  Given $\epsilon>0$, pick $n$ as in the second condition, and
  define a mapping $\Phi\colon\mathcal{F}\to \bbR^n$ by
  \[\Phi(x)=(x_1,\ldots,x_n).\]
  By the pointwise boundedness of $\mathcal{F}$, the image
  $\Phi(\mathcal{F})$ is totally bounded.

  If $x,y\in\mathcal{F}$ with
  $\abs{\Phi(x)-\Phi(y)}_p=\big(\sum_{k=1}^n\abs{x_k-y_k}^p\big)^{1/p}<\epsilon$, then
  \[
    \norm{x-y}_p\le
    \Big(\sum_{k=1}^n\abs{x_k-y_k}^p\Big)^{1/p}
   +\Big(\sum_{k>n}\abs{x_k-y_k}^p\Big)^{1/p}
   <\epsilon+2\epsilon=3\epsilon.
  \]
  By Lemma \ref{lemma:key}, $\mathcal{F}$ is totally bounded.

  \medskip\noindent
  We will leave proving the converse as an exercise to the reader.
  The techniques from the proof of Theorem \ref{thm:a} are easily
  adapted.  See also the proof of Theorem~\ref{thm:k-r}.
\end{proof}

\section{The Kolmogorov--Riesz theorem}

\begin{thm}[Kolmogorov--Riesz] \label{thm:k-r}
  Let $1\le p<\infty$.
  A subset $\mathcal{F}$ of $L^p(\bbR^n)$ is totally bounded if,
  and only if,
  \begin{enumerate}
  \item $\mathcal{F}$ is bounded,
  \item for every $\epsilon>0$ there is some $R$ so that, for
  every $f\in\mathcal{F}$,\[\int_{\abs{x}>R}\abs{f(x)}^p\,dx<\epsilon^p,\]
  \item for every $\epsilon>0$ there is some $\rho>0$ so that, for
  every $f\in\mathcal{F}$ and $y\in\bbR^n$ with $\abs{y}<\rho$,
 \[\int_{\bbR^n}\abs{f(x+y)-f(x)}^p\,dx<\epsilon^p.\]
  \end{enumerate}
\end{thm}
\begin{proof}
  Assume that $\mathcal{F}\subset L^p(\bbR^n)$ satisfies 
the three conditions.
  First, given $\epsilon>0$, pick $R$ as in the second condition,
  and $\rho$ as in the third condition.

  Let $Q$ be an open cube centered at the origin so that
  $\abs{y}<\frac12\rho$ whenever $y\in Q$.
  Let $Q_1,\ldots,Q_N$ be
  mutually non-overlapping translates of $Q$
  so that the closure of $\bigcup_i Q_i$ contains
  the ball with radius $R$ centered at the origin.
  Let $P$ be the projection map of $L^p(\bbR^n)$ onto the linear
  span of the characteristic functions of the cubes $Q_i$ given by
  \[
    Pf(x)=\begin{cases}
      \displaystyle\frac1{\abs{Q_i}}\int_{Q_i}f(z)\,dz,&
       x\in Q_i,\quad i=1,\ldots,N,\\
      0&\text{otherwise}.
    \end{cases}
  \]
  From \cond{ii} and the definition of $Pf$ we find, for $f\in\mathcal{F}$,
  \begin{align*}
    \norm{f-Pf}_p^p
    &<\epsilon^p+\sum_{i=1}^N\int_{Q_i}\abs{f(x)-Pf(x)}^p\,dx
  \\&=\epsilon^p+\sum_{i=1}^N\int_{Q_i}\abs[\Big]{\frac1{\abs{Q_i}}\int_{Q_i}
       \bigl(f(x)-f(z)\bigr)\,dz}^p\,dx.
  \end{align*}
  Next we use Jensen's inequality
    and change a variable of integration,
     where we note that $x-z\in2Q$ when $x,z\in Q_i$:
  \begin{align*}
    \norm{f-Pf}_p^p
    &<\epsilon^p+\sum_{i=1}^N\int_{Q_i}\frac1{\abs{Q_i}}\int_{Q_i}
       \abs[\big]{f(x)-f(z)}^p\,dz\,dx
  \\&\le\epsilon^p+\sum_{i=1}^N\int_{Q_i}\frac1{\abs{Q_i}}\int_{2Q}
       \abs[\big]{f(x)-f(x+y)}^p\,dy\,dx
  \\&\le\epsilon^p+\frac1{\abs{Q}}\int_{2Q}\int_{\bbR^n}
       \abs[\big]{f(x)-f(x+y)}^p\,dx\,dy
  \\&<\epsilon^p+\frac1{\abs{Q}}\int_{2Q}\epsilon^p\,dy
    =(2^n+1)\epsilon^p
  \end{align*}
  by \cond{iii}.  Thus $\norm{f-Pf}_p<(2^n+1)^{1/p}\epsilon$,
  and $\norm{f}_p<(2^n+1)^{1/p}\epsilon+\norm{Pf}_p$.
  By the linearity of $P$,
  if $f,g\in\mathcal{F}$ and $\norm{Pf-Pg}_p<\epsilon$ then
  $\norm{f-g}_p<\bigl((2^n+1)^{1/p}+1\bigr)\epsilon$.
  Moreover, since $P$ is bounded (in fact $\norm{P}=1$) and
  $\mathcal{F}$ is bounded by \cond{i},
  the image $P[\mathcal{F}]$ is bounded.
  Since the image of $P$ is finite dimensional,
  $P[\mathcal{F}]$ is totally bounded.
  Thus $\mathcal{F}$ is totally bounded by Lemma \ref{lemma:key}.

  \medskip\noindent
  For the converse, assume that $\mathcal{F}$ is totally bounded.

  The existence of a finite $\epsilon$-cover for $\mathcal{F}$,
  for any $\epsilon$, clearly implies the boundedness of
  $\mathcal{F}$, thus establishing Condition \cond{i}.

  To establish Condition \cond{ii}, let $\epsilon>0$ be given, let
  $\setof{U_1,\ldots,U_n}$ be an $\epsilon$-cover of $\mathcal{F}$,
  and chose $g_j\in U_j$ for $j=1,\ldots,n$.
  Select $R$ so that
  \[ \int_{x>R}\abs{g_j(x)}^p\,dx< \epsilon^p,\qquad j=1,\ldots,m. \]
  If $f\in U_j$ then $\norm{f-g_j}_p\le\epsilon$, and so
  \begin{align*}
     \Bigl(\int_{x>R}\abs{f(x)}^p\,dx\Bigr)^{1/p}
     &\le\Bigl(\int_{x>R}\abs{f(x)-g_j(x)}^p\,dx\Bigr)^{1/p}
        +\Bigl(\int_{x>R}\abs{g_j(x)}^p\,dx\Bigr)^{1/p}
   \\&\le\norm{f-g_j}_p+\Bigl(\int_{x>R}\abs{g_j(x)}^p\,dx\Bigr)^{1/p}
     <2\epsilon,
  \end{align*}
  thus establishing Condition \cond{ii}.

  Condition \cond{iii} is established similarly, by noting that the inequality
  of the condition is easily established for any \textit{single}
  function $f\in L^p(\bbR^n)$, for example using the fact that
  $C^\infty_c(\bbR^n)$ is dense in $L^p(\bbR^n)$.
  Then, picking an $\epsilon$-cover $\setof{U_1,\ldots,U_n}$
  and $g_j\in U_j$ for each $j$ as in the previous paragraph,
  given $\epsilon>0$ we can find $\rho>0$ with
  \[\int_{\bbR^n}\abs{g_j(x+y)-g_j(x)}^p\,dx<\epsilon^p,\qquad \abs{y}<\rho,\quad j=1,\ldots,m.\]
  Again, if $f\in U_j$ we find
  \begin{align*}
    \Bigl(\int_{\bbR^n}\abs{f(x+y)-f(x)}^p\,dx\Bigr)^{1/p}
    &\le\Bigl(\int_{\bbR^n}\abs{f(x+y)-g_j(x+y)}^p\,dx\Bigr)^{1/p}
    \\&\phantom{{}\le{}}+\Bigl(\int_{\bbR^n}\abs{g_j(x+y)-g_j(x)}^p\,dx\Bigr)^{1/p}
    \\&\phantom{{}\le{}}+\Bigl(\int_{\bbR^n}\abs{g_j(x)-f(x)}^p\,dx\Bigr)^{1/p}
  \\&<3\epsilon,
  \end{align*}
  and the proof is complete.
\end{proof}

\begin{remark}
  (I) A singleton set is clearly totally bounded,  yet Condition
  \cond{iii} is not obvious for a singleton set at first glance.
  However, it follows easily from the density of the space of smooth
  functions with compact support in $L^p$.

  (II) In applications, one sometimes constructs a sequence
  $f_1,f_2,\ldots$ in $L^p$ satisfying the first two conditions of
  Theorem \ref{thm:k-r} and the condition
  \[
    \Bigl(\int_{\bbR^n}\abs[\big]{f_n(x+y)-f_n(x)}\,dx\Bigr)^{1/p}
    <\alpha(y)+\beta(n),\qquad
    \lim_{y\to0}\alpha(y)=0,\quad\lim_{n\to\infty}\beta(n)=0.
  \]
  Then for some $N$ and $\delta>0$, the right-hand side of the
  above inequality is less than $\epsilon$ for all $n>N$ and
  $\abs{y}$ small enough.
  By the fact noted in the previous paragraph,
  we can choose a smaller upper bound for $\abs{y}$ to make the integral
  smaller than $\epsilon$ for $n=1,2,\ldots,N$.
  Thus $\setof{f_1,f_2,\ldots}$ satisfies Condition \cond{iii}, and
  hence a convergent subsequence exists.
\end{remark}

An interesting corollary to the Kolmogorov theorem is the following
result, see  \cite{Pego}, which also contains a variant using the uniform
smoothness of the functions in $\mathcal F$ and their Fourier transforms.
See also \cite{Dorfler}, which contains an alternate formulation based
on the short-time Fourier transform, as well as one based on the
wavelet transform.
\begin{cor} \label{coro}
Let $\mathcal{F}\subseteq L^2(\mathbb R^d)$ be such that
$\sup_{f\in \mathcal{F}}\norm{f}_2\le M<\infty$.
If
\[
\lim_{r\to\infty} \sup_{f\in \mathcal{F}}\int_{\abs{x}\ge r}\abs{f(x)}^2\, dx=0
\quad\text{and}\quad
\lim_{\rho\to\infty} \sup_{f\in \mathcal{F}}\int_{\abs{\xi}\ge \rho}\abs{\hat f(\xi)}^2\, d\xi=0,
\]
then $\mathcal{F}$ is totally bounded in $L^2(\mathbb R^d)$.
\end{cor}
\begin{proof}
  We show that $\mathcal{F}$ satisfies the conditions of Theorem
  \ref{thm:k-r} for $p=2$.
  Clearly, Conditions \cond{i} and \cond{ii} are among our
  assumptions, so we only need to prove \cond{iii}.
  For $f\in\mathcal{F}$ we find:
\begin{align*}
  \int_{\bbR^n}\abs{f(x+y)-f(x)}^2\,dx
    &=\int_{\bbR^n}\abs{(e^{i\xi\cdot y}-1)\hat f(\xi)}^2\,d\xi
  \\&\le\int_{\abs{\xi}<\rho}\abs[\big]{(e^{i\xi\cdot y}-1)\hat f(\xi)}^2\,d\xi
      +4\int_{\abs{\xi}>\rho}\abs[\big]{\hat f(\xi)}^2\,d\xi
   \\&\le M^2\sup_{\abs{\xi}<\rho}\abs{e^{i\xi\cdot y}-1}^2+\epsilon
   \qquad\text{for $\rho$ big enough}
   \\&<M^2\rho^2\abs{y}^2+\epsilon<2\epsilon
\end{align*}
  if $\abs{y}<\sqrt{\epsilon}/(\rho M)$.
  Here $\rho$, and hence the upper bound on $\abs{y}$,
  can be chosen independently of $f$.
  This shows Condition \cond{iii} of Theorem \ref{thm:k-r}
  and finishes the proof.
\end{proof}

In the following result,
$L^p\loc(\Omega)$ is equipped with the topology of $L^p$ convergence on
compact subsets of $\Omega$.
Recall that $\Omega$ is the countable union of compacts, e.g.,
$\Omega=K_1\cup K_2\cup\ldots$ with
$K_k=\setof{x\in\Omega\suchthat \abs{x}\le k\text{ and } \dist(x,\bbR^n\setminus\Omega)\ge1/k}$.
Moreover any compact subset of $\Omega$ is contained in some $K_k$,
and so the topology on $L^p\loc(\Omega)$ is given by the countable
family of seminorms $\norm{f}_k=\norm{f|_{K_k}}_{L^p(K_k)}$.
$L^p\loc(\Omega)$ is complete with respect to the metric
$(f,g)\mapsto\sum_{k=1}^\infty\min(2^{-k},\norm{f-g}_k)$.
\begin{cor}
  Let $\Omega\subseteq\bbR^n$ be an open set.
  Write $f_K(x)=f(x)$ when $x\in K$, $f_K(x)=0$ otherwise.
  A subset $\mathcal{F}\subseteq L^p\loc(\Omega)$ is totally bounded
  if, and only if, the following holds:
  \begin{enumerate}
  \item For every compact $K\subset\Omega$ there is some $M$ so that
  \[ \int\abs{f_K(x)}^p\,dx<M,\qquad f\in\mathcal{F}.\]
  \item For every $\epsilon>0$ and every compact $K\subset\Omega$
  there is some
  $\rho>0$ so that
  \[ \int\abs{f_K(x+y)-f_K(x)}^p\,dx<\epsilon^p,
  \qquad f\in\mathcal{F},\quad\abs{y}<\rho. \]
  \end{enumerate}
\end{cor}

\begin{proof}
  Note that $\mathcal{F}$ is totally bounded in $L^p\loc(\Omega)$
  if and only if
  $\mathcal{F}_k=\setof{f_{K_k}\colon f\in\mathcal{F}}$
  is totally bounded for every $k$,
  with $K_k$ as defined above.
\end{proof}

For the next result, recall that the \emph{Sobolev space} $W^{k,p}(\bbR^n)$
is defined to consist of those measurable functions $f$ which,
together with all their distributional derivatives
$D^\alpha f$ of order $\abs{\alpha}\le k$,
belong to $L^p(\bbR^n)$.
Here $\alpha=(\alpha_1,\dotsc,\alpha_n)$ is a multi-index,
i.e., each $\alpha_j$ is a nonnegative integer,
$\abs{\alpha}=\alpha_1+\dotsb+\alpha_n$,
and $D^\alpha=\partial^{\abs{\alpha}}/(\partial x_1^{\alpha_1}\cdots\partial x_n^{\alpha_n})$.
Finally, $W^{k,p}(\bbR^n)$ is equipped with the complete norm
\[\norm{f}_{k,p}=\Bigl(\int_{\bbR^n}\sum_{\abs{\alpha}\le k}
                    \abs{D^\alpha f(x)}^p\,dx\Bigr)^{1/p}.
\]

\begin{cor}
  A subset $\mathcal{F}\subseteq W^{k,p}(\bbR^n)$ is totally bounded
  if, and only if, the following holds:
  \begin{enumerate}
  \item $\mathcal{F}$ is bounded,
  i.e., there is some $M$ so that
  \[ \int \abs{D^\alpha f(x)}^p\,dx<M,\qquad f\in\mathcal{F},
     \quad\abs{\alpha}\le k.
  \]
  \item For every $\epsilon>0$ there is some $R$ so that
  \[\int_{\abs{x}>R}\abs{D^\alpha f(x)}^p\,dx<\epsilon^p,
    \qquad f\in\mathcal{F},\quad \abs{\alpha}\le k.\]
  \item For every $\epsilon>0$ there is some
  $\rho>0$ so that
  \[ \int_{\bbR^n}\abs{D^\alpha f(x+y)-D^\alpha f(x)}^p\,dx<\epsilon^p,
     \qquad f\in\mathcal{F},\quad\abs{\alpha}\le k,\quad\abs{y}<\rho.
  \]
  \end{enumerate}
\end{cor}

\begin{proof}
  Note that $\mathcal{F}$ is totally bounded in $W^{k,p}(\bbR^n)$
  if and only if
  $D^\alpha[\mathcal{F}]=\setof{D^\alpha f\colon f\in\mathcal{F}}$
  is totally bounded in $L^p(\bbR^n)$
  for every multi-index $\alpha$ with $\abs{\alpha}\le k$.
\end{proof}

\section{A bit of history}\label{history}

In 1931, Kolmogorov \cite{Kolmogorov} proved the first result in this
direction.
It characterizes compactness in $L^p(\bbR^n)$ for $1<p<\infty$, in
the case where all functions are supported in a common bounded set.
Condition \cond{iii} of Theorem \ref{thm:k-r} is replaced by the uniform
convergence in $L^p$ norm of spherical means of each function in the
class to the function itself.
(Clearly, our Condition \cond{ii} is automatic in this case.)

Just a year later, Tamarkin \cite{Tamarkin} expanded this result to
the case of unbounded supports by adding Condition \cond{ii} of Theorem
\ref{thm:k-r}.

In 1933, Tulajkov \cite{Tulajkov} expanded the Kolmogorov--Tamarkin
result to the case $p=1$.

In the same year, and probably independently, M.~Riesz \cite{Riesz}
proved the result for $1\le p<\infty$, essentially in the form of
our Theorem \ref{thm:k-r}.
Thus we feel somewhat justified in using the names Kolmogorov and
Riesz in referring to the theorem, though we are perhaps being a bit
unfair to Tamarkin and Tulajkov in doing so.

The compactness theorem has also seen generalizations in other
directions.

Hanson \cite{Hanson} proved a necessary and sufficient condition for
compactness of a family of measurable functions on a bounded
measurable set, with respect to convergence in measure.
(Here the measurable functions form a metric space in which the
distance between two functions is the infimum of all $\epsilon>0$ so
that the two functions differ by at most $\epsilon$ except on a set
of measure $\le\,\epsilon$.)

Fréchet \cite{Frechet} replaced Conditions \cond{i} and \cond{ii} of Theorem
\ref{thm:k-r} with a single condition (``equisummability''), and
generalized the theorem to arbitrary positive $p$.

Phillips \cite[Thm 3.7]{Phillips} proved a necessary and sufficient
condition for compactness in $L^p$ on a general measure space
($1\le p\le\infty$),
and indeed in \emph{any} Banach space,
which is however somewhat less suited to applications to PDEs.
Nevertheless, our sufficiency proof for Theorem \ref{thm:k-r} is based
on Phillips' criterion.
(It is more common, albeit more involved,
to use mollifiers in the proof.)

Weil \cite{Weil} (see also \cite[p.\ 269\thickspace ff]{Edwards}) extended the 
result to $L^p(G)$ where $G$ is a
locally compact group. Tsuji \cite{Tsuji} considered the case of $L^p(\bbR^d)$
with $0<p<1$, and Takahashi \cite{Takahashi} studied the same problem
in Orlicz spaces.  A characterization of compact subsets
of $L^p([0,T]; B)$ ($B$ a Banach space), which is very convenient in
the context of time-dependent partial differential equations, is given by Simon
\cite{Simon} (see also \cite{Maitre}).
A readable account of some of the historical development
can be found in \cite[p.\ 388]{DunfordSchwartz}. 
Helly's theorem \cite{Helly}, which was published already in
1912,  is easily seen to be a special case of
Kolmogorov's compactness theorem in the one-dimensional case, see
Section \ref{helly_sect}.

Further references include \cite{Veress},
\cite{Izumi}, \cite{BrunoGrande1}, \cite{BrunoGrande}, \cite{Feichtinger},
\cite{Nicolescu}.

\section{The Rellich--Kondrachov theorem}
In this section we use Kolmogorov's theorem to prove a simple variant
of the Rellich--Kondrachov theorem \cite{Rellich, Kondrachov}.
Our simplification consists in avoiding boundary regularity conditions
by working on the entire space $\bbR^n$.
The standard Rellich--Kondrachov theorem requires a bounded region.
The present version replaces this by a uniform decay estimate,
specially tailored to fit the framework of the present paper.

The Sobolev norm $\norm{f}_{1,p}$ on $W^{1,p}(\bbR^n)$ is defined by
\[
  \norm{f}_{1,p}
  =\Bigl(\int_{\bbR^n}\bigl(\abs{f(x)}^p+\abs{\nabla f(x)}_p^p\bigr)\,dx
   \Bigr)^{1/p},\qquad
  \abs{\nabla f}_p=\Bigl(\sum_{j}
   \abs[\Big]{\frac{\partial f}{\partial x_i}}^p\Bigr)^{1/p}.
\]
According to the Sobolev embedding theorem, if $p<n$ then
$W^{1,p}(\bbR^n)\subset L^q(\bbR^n)$, and the inclusion map is bounded,
for any $q$ satisfying $p\le q\le p^*$,
where $p^*$ is the \emph{conjugate Sobolev exponent}:
\[
  \frac1{p^*}=\frac1p-\frac1n.
\]
To see where this exponent comes from, consider a function $f$ and its scalings
$f^\lambda(x)=f(x/\lambda)$ where $\lambda>0$, and note that
$\norm{f^\lambda}_p=\lambda^{n/p}\norm{f}_p$ and
$\norm{\nabla f^\lambda}_p=\lambda^{n/p-1}\norm{\nabla f}_p$,
so the inclusion map $W^{1,p}\to L^q$ can only be bounded if
there exists a constant $C$ with
$\lambda^{n/q}\le C(\lambda^{n/p}+\lambda^{n/p-1})$
for all $\lambda>0$.
In the limits $\lambda\to\infty$ and $\lambda\to0$
we conclude $n/q\le n/p$ and $n/q\ge n/p-1$ respectively.

\begin{thm}
Assume $p<n$ and $p\le q<p^*$,
and let $\mathcal{F}$ be a bounded subset of $W^{1,p}(\bbR^n)$.
Assume that for every $\epsilon>0$ there exists some $R$ so that, for every
$f\in\mathcal{F}$,
\[
  \int_{\abs{x}>R}\bigl(\abs{f(x)}^p+\abs{\nabla f(x)}_p^p\bigr)\,dx<\epsilon^p.
\]
Then $\mathcal{F}$ is a totally bounded subset of $L^q(\bbR^n)$.
\end{thm}

\begin{proof}
We shall show that $\mathcal{F}$ satisfies the hypotheses of
Theorem \ref{thm:k-r}, with $p$ replaced by $q$.
We shall use the Sobolev embedding inequality
$\norm{f}_q\le C\norm{f}_{1,p}$,
where the constant $C$ depends only on $p$, $q$ and $n$,
and which is valid under the stated assumption,
see \cite[ 4.30 (p. 101) and Theorem 4.12 I C (p. 85) with $j=0$,
  $k=n$, $m=1$]{AdamsFournier}.
Condition \cond{i} of Theorem \ref{thm:k-r} follows immediately
from the Sobolev embedding inequality.
Condition \cond{ii} is almost equally immediate,
from applying the Sobolev embedding inequality to the function
$x\mapsto f(x)\chi(\abs{x}-R)$, where $\chi\in C^\infty(\bbR)$,
$0\le\chi\le 1$, $\chi(x)=0$ for $x<0$ and $\chi(x)=1 $ for $x>1$.

If we apply the Sobolev embedding inequality to the function
$x\mapsto f(x/\lambda)$ where $\lambda>0$ and change variables in the
resulting integrals, we obtain
\begin{equation}\label{eq:sobolevrescaled}
  \lambda^{n/q}\norm{f}_q
  \le C
  \Bigl(\lambda^n\int_{\bbR^n}\abs{f(x)}^p\,dx
        +\lambda^{n-p}\int_{\bbR^n}\abs{\nabla f(x)}_p^p\,dx
   \Bigr)^{1/p}
\end{equation}
We shall apply the above inequality not to $f$,
but to $x\mapsto f(x+y)-f(x)$, where $f\in\mathcal F$.

Now let $\epsilon>0$ be given.
By picking $\lambda$ sufficiently large we can ensure that
\begin{equation}\label{eq:sobdiff}
  C\Bigl(\lambda^{n-p}\int_{\bbR^n}\abs{\nabla f(x+y)-\nabla f(x))}_p^p\,dx\Bigr)^{1/p}
  \le\epsilon\lambda^{n/q}
\end{equation}
for all $f\in\mathcal F$,
since the integral in this expression is bounded
uniformly for $f\in\mathcal F$.

Next, we find (using the Jensen and Hölder inequalities, then Fubini's
theorem)
\begin{align*}
  \int_{\bbR^n}\abs{f(x+y)-f(x)}^p\,dx
  &=\int_{\bbR^n}\abs[\Big]{\int_0^1 y\cdot\nabla f(x+ty)\,dt}^p\,dx
\\&\le\abs{y}_{p'}^p\int_0^1\int_{\bbR^n}\abs{\nabla f(x+ty)}_p^p\,dx\,dt
\\&=\abs{y}_{p'}^p\int_{\bbR^n}\abs{\nabla f(x)}_p^p\,dx
\end{align*}
(where $p$ and $p'$ are conjugate exponents)
for any test function $f$, and hence for any $f\in W^{1,p}$.
The integrals on the right-hand side of this inequality are uniformly
bounded for $f\in\mathcal{F}$,
and so we can find some $\delta>0$ so that $\abs{y}<\delta$
implies
\begin{equation}\label{eq:sobfunc}
  C\Bigl(\lambda^{n}\int_{\bbR^n}\abs{f(x+y)-f(x)}^p\,dx\Bigr)^{1/p}
  \le\epsilon\lambda^{n/q}.
\end{equation}
For such $y$ and $f$, \eqref{eq:sobolevrescaled} applied to
$x\mapsto f(x+y)-f(x)$
combined with \eqref{eq:sobdiff} and \eqref{eq:sobfunc} to yield
\[
  \lambda^{n/q}\norm{f(\cdot+y)-f(\cdot)}_q\le2^{1/p}\epsilon\lambda^{n/q},
\]
and so assumption (\textit{iii}) of Theorem \ref{thm:k-r} is satisfied.
\end{proof}

\section{Helly's theorem}\label{helly_sect}

Helly's theorem is often referred to as
\emph{Helly's selection principle},
in order to avoid confusion with another theorem by Helly,
stating that, given a collection of convex sets in $\bbR^n$
so that any $n+1$ of them have a point in common,
then any finite subcollection has nonempty intersection.
Helly's selection principle is essentially a corollary of the Kolmogorov--Riesz
theorem, though historically it was not derived that way.

Recall that an integrable function $f$ on the line
is of bounded variation if it has finite essential or total variation,
that is, if
\[
  \TV(f) =\sup \sum_{j=1}^{m}\abs{f(x_{j+1})-f(x_j)}<\infty,
\]
where the supremum is taken over all finite partitions $x_j<x_{j+1}$
such that each $x_j$ is a point of approximate continuity of $f$
(that is, $\delta^{-1}\abs{\setof{x\suchthat\abs{x-x_j}<\delta,
 \abs{f(x)-f(x_j)}\ge\epsilon}}\to0$
for every $\epsilon>0$ as $\delta\to0$. See, e.g., \cite[p.~47]{Evans}).
We need a lemma:
\begin{lem} \label{bv-trans}
  Let $u$ be function of bounded variation on $\bbR$.  Then
\[
  \int_{-\infty}^\infty\abs{u(x+y)-u(x)}\,dx\le\abs{y}\TV(u)
\]
  for all $y\in\bbR$.
\end{lem}
\begin{proof}
We may assume $y>0$ without loss of generality.
The calculation
\begin{align*}
  \int_{-\infty}^\infty\abs{u(x+y)-u(x)}\,dx
   &=\sum_{j=-\infty}^\infty\int_0^y \abs{u(x+jy+y)-u(x+jy)}\,dx
 \\&=\int_0^y\sum_{j=-\infty}^\infty \abs{u(x+(j+1)y)-u(x+jy)}\,dx
 \\&\le\int_0^y\TV(u)\,dx=y\TV(u),
\end{align*}
finishes the proof.
\end{proof}
\begin{thm}[Helly]
  Let $(u_n)$ be a sequence of functions of bounded variation on the
  bounded real interval $[a,b]$.  If there is a constant $M$ so that
  $\TV(u_n)\le M$ and $\norm{u_n}_\infty\le M$ for all $n$, then there is a
  subsequence of $(u_n)$ which converges pointwise everywhere
  and in $L^1$ norm in $[a,b]$ to a
  function of bounded variation.
\end{thm}
\begin{proof}
Extend each function $u_n$ to all of $\bbR$ by setting it to zero
outside $[a,b]$.  By Lemma \ref{bv-trans}, the set of all the
functions $u_n$ satisfy Condition \cond{iii} of Theorem \ref{thm:k-r} (with
$p=1$), while $(i)$ holds by assumption and $(ii)$ is trivial.  Hence
there is a subsequence of $(u_n)$ which converges in $L^1([a,b])$.
Moreover, integration theory tells us
that we also get pointwise convergence almost everywhere, possibly
after passing to a subsequence once more.
However, this is not quite enough.

Write instead $u_n=v_n-w_n$ where each $v_n$, $w_n$ is an
non-decreasing function: $v_n(x)$ is $u_n(a)$ plus the positive
variation of $u_n$ on
the interval $[a,x]$, and  $w_n(x)$ is the negative variation on the
same interval.  Then the sequences $(v_n)$ and $(w_n)$ both satisfy
the conditions of the present theorem, and so, by the result of the
previous paragraph, we may pass to a subsequence so that $(v_n)$ and
$(w_n)$ both converge in $L^1([a,b])$,
as well as pointwise almost everywhere.

Let $v$ be the limit of the sequence $(v_n)$.  Clearly, $v$ is
non-decreasing on the set where pointwise convergence holds,
and so we may assume that $v$ is non-decreasing everywhere, after
possibly redefining it on a set of measure zero.

Now it is clear that $v_n(x)\to v(x)$ for any point of continuity $x$
for $v$: Given $\epsilon>0$, pick $\delta>0$ so that
$\abs{y-x}<\delta$ implies $\abs{v(y)-v(x)}<\epsilon$,
let $x-\delta<y<x<z<x+\delta$ with $v_n(y)\to v(y)$ and $v_n(z)\to
v(z)$, and note that for $n$ large enough we get
$v(x)-2\epsilon<v(y)-\epsilon<v_n(y)
 \le v_n(x)
 \le v_n(z)<v(z)+\epsilon<v(x)+2\epsilon$,
so that $\abs{v_n(x)-v(x)}<2\epsilon$.

Since $v$ has at most a countable number of discontinuities,
a diagonal argument yields a further subsequence which converges at
all the discontinuities of $v$ as well, and so we have pointwise
convergence everywhere.

In the same way we show that $w_n(x)\to w(x)$ for all $x$.
Thus $u_n\to v-w$ pointwise, and $v-w$ has bounded variation.
\end{proof}

\begin{remark}
The above proof is probably not the most natural one,
but it does make clear the connection with the Kolmogorov--Riesz theorem.
In a sense $L^1$ convergence is irrelevant:
Pointwise convergence is the key,
and $L^1$ convergence follows from the bounded convergence theorem.

It should be noted, however, that Helly's theorem,
without pointwise convergence,
is also true in higher dimensions \cite[p.~176]{Evans}.

A recent generalization of Helly's selection principle
(in one dimension) can be found in
\cite{TretyachenkoChistyakov}.
\end{remark}


\end{document}